\documentclass[11pt]{amsart}
\usepackage{amsmath,amssymb, amsthm}

\newcommand\Z{{\mathbb Z}}
\newcommand\R{{\mathbb R}}

\newcommand\C{{\mathbb C}}
\newcommand\N{{\mathbb N}}
\newcommand\T{{\mathbb T}}
\newcommand\g{\gamma}
\newcommand\G{\Gamma}
\newcommand\inv{^{-1}}
\newtheorem{theorem}{Theorem}[section]

\newtheorem{lemma}[theorem]{Lemma}
\newtheorem{corollary}[theorem]{Corollary}
\newtheorem{defn}[theorem]{Definition}
\newtheorem{example}[theorem]{Example}
\newtheorem{question}[theorem]{Question}

\DeclareMathOperator{\Fix}{Fix} \DeclareMathOperator{\Out}{Out}
\DeclareMathOperator{\Int}{Int}  \DeclareMathOperator{\Aut}{Aut}
\DeclareMathOperator{\Diff}{Diff}\DeclareMathOperator{\Inn}{Inn}
\DeclareMathOperator{\PInn}{PInn}

\begin{document}

\title{When is a group action determined by it's orbit structure?}
\author{David Fisher and Kevin Whyte}
\thanks{First author partially supported by NSF grants DMS-9902411 and DMS-0226121.
Second author partially supported by NSF grants DMS-9971094 and
DMS-0204576.}
\begin{abstract}
We present a simple approach to questions of topological orbit
equivalence for actions of countable groups on topological and
smooth manifolds. For example, for any action of a countable group
$\Gamma$ on a topological manifold where the fixed sets for any
element are contained in codimension two submanifolds, every orbit
equivalence is equivariant.  Even in the presence of larger fixed
sets, for actions preserving rigid geometric structures our
results force sufficiently smooth orbit equivalences to be
equivariant. For instance, if a countable group $\Gamma$ acts on
${\mathbb T}^n$ and the action is $C^1$ orbit equivalent to the
standard action of $SL_n({\mathbb Z})$ on ${\mathbb T}^n$, then
$\Gamma$ is isomorphic to $SL_n({\mathbb Z})$ and the actions are
isomorphic. (The same result holds if we replace $SL_n(\mathbb Z)$
by a finite index subgroup.) We also show that preserving a
geometric structure is an invariant of smooth orbit equivalence
and give an application of our ideas to the theory of hyperbolic
groups.

In the course of proving our theorems, we generalize a theorem of
Sierpinski which says that a connected Hausdorff compact
topological space is not the disjoint union of countably many
closed sets.  We prove a stronger statement that allows "small"
intersections provided the space is locally connected.  This
implies that for any continuous action of a countable group $\G$
on a connected, locally connected, locally compact, Hausdorff
topological space, where the fixed set of every element is
"small", every orbit equivalence is equivariant.
\end{abstract}

\maketitle

\section*{Introduction}

   Many interesting properties of a dynamical system are properties of
the orbits: minimality, periodic points, etc.  In this paper we
study the inverse problem; to what extent does the orbit structure
determine the action?

In the measurable category, the topic of orbit equivalence of
group actions has been studied in great detail.  For amenable
groups, measurable orbit equivalence of ergodic measure preserving
actions is shown to be a trivial relation in \cite{CFW}.  For
non-amenable groups, many interesting rigidity theorems have been
proven, starting with Zimmer, and continuing with work of Adams,
Furman, Gaboriau, and Monod-Shalom \cite{Ad, Fu,G, MS,Z}. All of
these results depend heavily on deep machinery, whether it is
Zimmer's cocycle superrigidity and Ratner's measure classification
theorem or the theories of $L^2$ or bounded cohomology.  In this
note, we prove analogues and variants of results known in the
measurable category in the topological and smooth categories. The
most striking fact is that here we can prove analogous theorems by
almost entirely elementary methods and without assuming an
invariant measure or volume for the action.

\begin{defn}
\label{definition:orbitequivalent} Let $X$ and $X'$ be topological
spaces and $\Gamma$ and $\Gamma'$ groups.  Two actions
$(X,\Gamma)$ and $(X',\Gamma')$ are {\em topologically orbit
equivalent} if there is a homeomorphism $f$ from $X$ to $X'$ which
maps the orbit relation for $\Gamma$ to the orbit relation for
$\Gamma'$.  The map $f$ is called a {\em topological orbit
equivalence}.
\end{defn}

Since we only consider topological orbit equivalence in this
paper, we will use the phrase orbit equivalent and orbit
equivalence to mean topological orbit equivalent and topological
orbit equivalence.

\begin{defn}
\label{definition:orbitrigid} Let a group $\Gamma$ act on a
topological space $X$.  We say that the action is {\em $C^0$ OE
rigid} if any orbit equivalence from $(X',\Gamma')$ to
$(X,\Gamma)$ is equivariant.

Furthermore, if $X$ is a smooth manifold and $\Gamma$ acts
smoothly, we say the action is {\em $C^k$ OE rigid} if any orbit
equivalence which is a $C^k$ diffeomorphism is equivariant.
\end{defn}

Given a group $\Gamma$ acting on a space $X$ and
$\gamma{\in}\Gamma$, we let $\Fix(\gamma)$ be the set of $\gamma$
fixed points.

\begin{theorem}
\label{theorem:codim2} Let $\Gamma$ be a countable group acting on
a compact manifold $X$.  Assume that $\Fix(\gamma)$ is contained
in a submanifold of codimension two for every $\gamma{\in}\Gamma$.
Then the $\Gamma$ action on $X$ is $C^0$ OE rigid.
\end{theorem}

Both free actions and complex analytic actions trivially satisfy
the hypotheses of the theorem.  See section
\ref{section:smallfixedsets} for more examples.  After proving
this result, we discovered that some related prior related results
for $\mathbb Z$ actions see \cite[Remark 3.4]{BT} and
\cite{K,GPS}.  All of these results follow from Theorem
\ref{small} below.   This theorem is proved using the
generalization of Sierpinski's theorem mentioned in the abstract.

Even in cases where fixed sets are larger, for example codimension
one submanifolds, in the presence of a rigid geometric structure
we can show that any non-equivariant orbit equivalence is not
(very) smooth.  For example:

\begin{theorem}
\label{theorem:slnz1} Let $\Gamma<SL_n(\mathbb Z)$ be a subgroup
of finite index.  Then the standard $\Gamma$ action on ${\mathbb
T}^n$ is $C^1$ OE rigid.
\end{theorem}

In this setting our techniques, combined with a few simple
observations, can easily identify all self orbit equivalences of
the action of $\Gamma$ on ${\mathbb T}^n$.

In addition, we see that preserving a geometric structure is an
invariant of sufficiently smooth orbit equivalences.  See section
\ref{section:geometricstructures} for definitions and discussion.

After essentially completing the work presented here, we
discovered the large body of work on orbit equivalence of $\mathbb
Z$ actions and the associated $C^*$ algebras, see for example
\cite{BH,BT,GPS} as well as the references mentioned above.  It
would be interesting to study the analogous questions regarding
associated group transformation $C^*$ algebras and orbit
equivalences for groups larger than $\mathbb Z$.

This project began as a conversation between the authors and R.
Spatzier at the Newton Institute of Mathematical Sciences. We
thank Spatzier for stimulating conversations and the Institute for
it's support.

\section{Basic Techniques}
\label{section:basictecniques}
 Given a space $X$ and a group $\G$
acting on $X$, we define $\Aut(X,\G)$ to be the group of
homeomorphisms of $X$ preserving the orbit equivalence relation,
and $\Inn(X,\G)$ to be the subgroup of homeomorphisms which send
every $x$ to a point in $\G x$. Clearly $\G \subset \Inn(X,\G)$.
When $X$ is a smooth manifold, $\Aut^k(X,\G)$ and $\Inn^k(X,\G)$
refer to the subgroups of $C^k$-diffeomorphisms in $\Aut(X,\G)$
and $\Inn(X,\G)$.  Below we abuse notation by letting
$\Inn^0(X,\G)=\Inn(X,\G)$ and $\Aut^0(X,\G)=\Aut(X,\G)$ whether
$X$ is smooth or not.

The prototypical result of the sort we are after is that
$\G=\Inn^k(X,\G)$.  When this holds, it implies that for $(Y,\G')$
orbit equivalent to $(X,\G)$,  there is an isomorphism $\G'  \to
\G$ which makes the given orbit equivalence equivariant. In
particular, it follows that $\Aut^k(X,\G)=N_{\Diff^k(X)}(\G)$, the
normalizer of $\Gamma$ in $\Diff^k(X)$. Actions for which this
holds are determined by their orbit structure in the strongest
possible sense.

We note that there is a short exact sequence:
$$1{\rightarrow}Z_{\Diff^k(X)}(\G){\rightarrow}N_{\Diff^k(X)}(\G){\rightarrow}\Aut(\G){\cap}\Diff^k(X){\rightarrow}1$$
where $Z_{\Diff^k(X)}(\G)$ denotes the centralizer in $\Diff^k(X)$
of $\Gamma$ and $\Aut(\G){\cap}\Diff^k(X)$ denotes those
automorphisms of $\G$ which can be realized as restrictions of
inner automorphisms of $\Diff^k(X)$.  It is frequently possible to
identify both $Z_{\Diff^k(X)}(\G)$ and $\Aut(\G){\cap}\Diff^k(X)$
explicitly and so obtain an exact description of $\Aut^k(X,\G)$.

Given $X$ and $\G$, and an inner orbit equivalence $f$, define the
sets
$$T_\g =\{x | f(x)=\g x\}$$

Clearly these sets are closed, and by definition they cover $X$.
The intersection of $T_\g$ and $T_{\g'}$ is contained in the fixed
set of $\g^{-1}\g'$.  On each $T_\g$, $f$ agrees with translation
by $\g$, so $f$ is, in this sense "piecewise $\G$".  Conversely,
anything which is piecewise in $\G$ is an inner orbit equivalence.
The following example was shown to us by Scot Adams.  A very
similar example occurs as \cite[Example 4.3]{BT}.

\begin{example}
\label{adams}
Let $X=\T^2=\R^2/\Z^2$ and $\G=SL_2(\Z)$. Let $f$ be the identity
for $0 \leq x \leq \frac{1}{2}$, and let $f(x,y)=(x,y+2x)$ for
$\frac{1}{2}\leq x \leq 1$.  It is easy to check that this is a
homeomorphism of the torus, and it is clearly an inner orbit equivalence.

\end{example}

This example is piecewise $SL_2(\Z)$ in a much nicer sense than
the above.  For a general $X$ and $\G$ we define the piecewise
$\G$ maps, $\PInn(X,\G)$ as those inner orbit equivalences for
which a finite collection of $T_{\g}$ cover $X$.  When the orbit
structure is not rigid, one can still hope that
$\Inn(X,\G)=\PInn(X,\G)$ as the later are somewhat more easily
understood.

Let $X=S^1 = \R \cup \{\infty\}$, and let $\G$ be the group of
rational affine maps.  Let $f$ be the map which is the identity
outside of $(0,2)$, sends $\frac{1}{n}$ to $\frac{1}{n+1}$ for all
$n \in \N$, and which is affine on the intervals
$[\frac{1}{n+1},\frac{1}{n}]$.  This is a (Lipschitz) inner orbit
equivalence which is not in $\PInn(X,\G)$.  This is approximately
as bad as it gets for countable groups.

\begin{lemma}
\label{lemma:key} If $X$ is a complete metric space and $\G$ is
countable, then $\cup_\G Int(T_\g)$ is an open dense subset of
$X$.
\end{lemma}

\begin{proof}
This is essentially the Baire category theorem.  Let $U$ be any non-empty
open set in $X$.  If $U$ is disjoint from $\cup_\G Int(T_\g)$ then
$U=\cup_G \partial (T_\g)$, which is a countable union of closed, nowhere
dense sets.  This is impossible by the Baire category theorem, hence
$U$ intersects $\cup_\G Int(T_\g)$, which proves density.  Openness is
obvious.
\end{proof}

Thus, on an open dense set, every inner orbit equivalence behaves
locally like an element of $\G$.  This is surprisingly powerful in
many cases.

\section{Geometric Structures}
\label{section:geometricstructures}

Many of the most natural examples of group actions preserve some
extra structure.   One unified treatment of these extra structures
is via "geometric structures a la Gromov", which are a
generalization of a reduction of the structure group of the frame
bundle.  Let $D^k$ be the group of $k$-jets of diffeomorphisms of
$\R^n$ preserving the origin, and $P^k(M)$ the $k^{th}$ order
frame bundle of $M$, viewed as a principle $D^k$ bundle.

\begin{defn}
\label{defn:gs} A {\em geometric structure} on a manifold $M$
consists of the following data:
\begin{enumerate}

\item a manifold $V$

\item a $D^k$ action on $V$

\item a $D^k$ equivariant map $\omega:P^k(M){\rightarrow}V$.
\end{enumerate}
We frequently refer to $\omega$ as the geometric structure leaving
$V$ and $k$ implicit. We say two geometric structures are of the
{\em same type} if the $k$'s are the same and the $V$'s are
isomorphic as $D^k$ spaces.
\end{defn}

Basic examples are familiar: metrics, volume forms, symplectic
structures, affine connections, projective connections, etc. The
reader unhappy with the above definition will lose little by
considering only the above list.  Note that we are not assuming
that our geometric structures are $A$-structures in the sense of
Gromov, as is necessary for most applications.

\begin{theorem}[Orbit Equivalences Preserve geometric Structures]
\label{oegs} Let  $M$  be a compact manifold with a $C^k$-action
of a countable group $\G$, preserving a geometric structure
$\omega$ of order $k$. Then $\Inn^k(M,\G) \subset \Aut^k(\omega)$.
\end{theorem}

\begin{proof}
Let $f$ be the orbit equivalence of $M$.  Let $\omega'$ be the pullback of
$\omega$, in other words the composition of the map $f^k: P^k(M) \to P^k(M)$ with
the map $P^k(M) \to V$ defining $\omega$.   By the key lemma,  $f$ agrees locally
on an open dense set with elements of $\G$.  Since $\G$ preserves $\omega$, this
means that $\omega$ and $\omega'$ are equal on an open dense set.  By continuity,
they are equal.
\end{proof}

\begin{corollary}
Let $M$ be a compact manifold with $\G$ action preserving a $k^{th}$-order geometric
structure.  Let $(N,\G')$ be any manifold with an action $C^k$ orbit equivalent to $(M,\G)$.
The $\G'$ action on $N$ preserves a geometric structure of the same kind as $\omega$.
\end{corollary}

\begin{proof}
The rigid structure on $N$ preserved by $\G'$ is the pullback of
$\omega$ via the $C^k$ orbit equivalence $f:N{\rightarrow}M$.
Since, for $\gamma'{\in}\G'$, the diffeomorphism
$f{\circ}\gamma'{\circ}f{\inv}$ of $M$ is in $\Inn^k(M,\G)$, the
result follows.
\end{proof}
This implies, for example, that one can recognize whether an action is volume preserving
or symplectic just from the orbit structure.  In some cases one can see how to do this, for
example, one can try to build a measure by counting to fraction of periodic orbits in a set.
If the finite orbits are equidistributed, this will recover the volume.  Of course this does not
always work.   It is unclear how to even begin to build a symplectic form out of the orbits,
although the Arnold conjectures give a way to recognize that certain orbit structures cannot
arise from symplectic actions.

When the automorphism group of a geometric structure is smaller,
for example for isometric or affine structures, Theorem \ref{oegs}
is very strong.   One class of such structures are the ones Gromov
calls rigid.  For our purposes the most convenient formulation of
this is:

\begin{defn}
\label{definition:rgs} If $\omega$ is a geometric structure on
$M$, $\omega$ is called {\em rigid} (of order $r$) if the action
of $\Aut^r(\omega)$ is free and proper on $P^r(M)$.
\end{defn}

This notion of a rigid geometric structure is equivalent to the
one defined in \cite{Gr1}, as shown there on pages 69-70.

\begin{theorem}[Actions Preserving Rigid Structures Are Rigid]
\label{theorem:rigidity} Let $\Gamma$ be a countable group acting
on a compact manifold $M$.  Assume the action preserves a rigid
geometric structure, $\omega$, of order $r$.   Then $\G =
\Inn^r(M,\G)$.
\end{theorem}

\begin{proof}
Since $\omega$ is rigid, $\Aut^r(\omega)$ is a Lie group which
acts on $M$, freely and properly on $P^r(M)$.  By theorem
\ref{oegs}, any inner  orbit equivalence is an element of
$\Aut^r(\omega)$. Suppose $g \in \Aut^r(\omega)$ preserves $\G$
orbits.   Then $M = \cup_\G T_\g$, with $gm = \g m$ for all $m \in
T_\g$.  If $T_\g$ has non-empty interior then there is some open
set on which $g$ and $\g$ agree, in particular, they agree at some
point of $P^r(M)$. Since the $\Aut^r(\omega)$ action is free on
$P^r(M)$, this implies $\g=g$.
\end{proof}

\begin{corollary}
Let $\Gamma$ be a countable group acting on a  compact manifold
$M$.  Assume the action preserves a rigid geometric structure,
$\omega$, of order $r$.   Then $\Aut^r(M,\G)$ consists only of
equivariant maps, and therefore $\Aut^r(M,\G)=N_{\Diff^r(X)}(\G)$.
\end{corollary}

Recall that $N_{\Diff^r(X)}(\G)$ is an extension of
$Z_{\Diff^r(X)}(\G)$ by a subgroup of $\Aut(\G)$.  For many $\G$
$\Out(\G)$ is known to be finite, for example all irreducible
lattices in semi-simple Lie groups not locally isomorphic to
$SL_2(\mathbb R)$, and all one ended hyperbolic groups. In
addition, for sufficiently hyperbolic actions the centralizer can
often be shown to be trivial.  For such groups the corollary
implies $\G$ is finite index in $\Aut^r(M,\G)$.  In specific cases
one can often identify the group $\Aut^r(M,\G)$ completely.

\begin{theorem}
\label{theorem:slnz} Let $\G<SL_n(\mathbb Z)$ be of finite index.
For the standard action of $\G$ on ${\mathbb T}^n$, the group
$\Aut^1({\mathbb T}^n, \G)$ is of finite index in $GL_n(\mathbb
Z)$.
\end{theorem}

For $n{\geq}3$ one can use Mostow rigidity to calculate
$\Aut(SL_n(\Z))$, but by considering the action on $H^1(M,{\mathbb
Z})$ it is easy to see that no automorphism not in $GL_n(\mathbb
Z)$ can be realized as a diffeomorphism of the torus. In
particular, the theorem holds even for $n=2$, despite the fact
that $\Aut(SL_2(\Z))$ is quite large. Easy dynamical arguments
using the presence of hyperbolic matrices in $\G$ show that the
centralizer of the action is trivial.

Note that Theorem \ref{theorem:slnz1} of the introduction follows
easily from this one.

One can state much more general results concerning affine actions
of lattices on homogeneous spaces.  Here one can use Ratner's
theorem to identify centralizers and Mostow rigidity to identify
automorphisms. If the action is sufficiently hyperbolic, the use
of Ratner's theorem can be replaced with techniques deriving from
\cite{HPS}. In particular, most of the results of \cite{Fu2} can
be reproven in the $C^1$ (rather than measurable) category by
these methods.

The differentiability hypotheses in Theorem \ref{theorem:rigidity}
are necessary in general. The group of $C^0$ orbit equivalences of
$SL_n(\Z)$ action on the torus is large. See section \ref{groups}
for some discussion of this case.   This is not simply a case of
continuous versus differentiable, as the higher order derivatives
are also needed:
\begin{example}
\label{example:circle} Let $X=S^1=\R \cup \{\infty\}$ and
$\G=SL_2(\Z[\frac{1}{2}])$. Let $f$ be the map which is the
identity for $x<0$, $f(x)=\frac{2x}{2-x}$ on $[0,1]$,
$6-\frac{4}{x}$ on $[1,2]$, and $x+2$ for $x>2$.  It is easy to
check that this is a $C^1$-inner orbit equivalence.  Since the
action of $\G$ preserves a rigid structure of order $2$, one sees
that the differentiability assumptions are needed.  This example
can be suspended to give similar examples on $\R P^n$ for all $n$.
\end{example}

\section{Small Fixed Sets}
\label{section:smallfixedsets}

% \marginpar{Do we need Hausdorff or
%something stronger?   The exact separation we use repeatedly is:
%given a point and a closed set, there is an open set around the
%point which is disjoint from the closed set. This seems to follow
%from Hausdorff if the closed set is compact, but not otherwise.}

The last section explored geometric reasons for an orbit structure
to be rigid.  There are also topological/dynamical reasons.    The
prototype is the case of free actions. If $\G$ acts freely on $M$,
then the decomposition $M = \cup_\G T_\g$ is a decomposition into
a countable disjoint closed sets.  If $M$ is connected, compact
and Hausdorff this implies, by a theorem of Sierpinski \cite{Si},
that all but one of these sets is empty.  Thus, for any free
action the orbit structure is rigid.  We generalize that result to
cases where the fixed sets are small. It is interesting to note
that in the context of the previous section, and of lattices in
Lie groups, the elements with big fixed sets are unipotents. Thus,
in this case, the presence of unipotents in a subgroup is an
obstruction to rigidity.

Throughout this section, we let the topological space $X$ be
connected, locally connected, Hausdorff, and locally compact.

A subset $A$ in $X$ {\em locally disconnects} if there is a connected open
set $U$ such that $U \setminus A$ is not connected.

A subset $Z$ of $X$ will be called {\em small} if it is closed,
has empty interior, and does not locally disconnect.

\begin{itemize}

\item Submanifolds of codimension 2 and higher are small.
\item $X$ has no cut points iff points are small.

\end{itemize}

\begin{lemma}
\label{lemma:trivial}
 A finite union of small sets in small.
\end{lemma}

\begin{proof}
If $A$ and $B$ are small, then clearly $C=A \cup B$ is closed with empty
interior.  Suppose $C$ disconnects a connected open set $U$, so that
$U \setminus C = V \cup W$, disjoint open sets.  If $B \subset
A$ this contradicts $A$ being small.  Let $b$ be a point of $B \setminus A$,
and let $N$ be a tiny connected open neighborhood of $b$ which is disjoint
from $A$.  Since $B$ is small, it cannot disconnect $N$, so $N \setminus B$
is contained in a single component of $U \setminus C$.  Let $\hat{V}$
(resp. $\hat{W}$) be $V$ (resp. $W$) union those points of $B \setminus A$
whose neighborhoods only intersect $V$ (resp. $W$).  These are disjoint
open sets whose union is $U \setminus A$.  This contradicts the fact that
$A$ is small.
\end{proof}

\begin{theorem} Let $X$ be connected, locally connected, Hausdorff and locally compact.
Let $Y_i$ be a countable collection of closed sets which cover $X$.  If each
intersection $Y_i \cap Y_j$ is small, then $Y_i = X$ for some $i$.
\end{theorem}

\begin{proof}

Without loss of generality, we may assume that $Int(Y_i) \cap Y_j$
is empty whenever $i \neq j$ (just succesively replace $Y_i$ by
$Y_i \setminus (Y_i {\cap} \cup_{j > i} Int(Y_j)$ for
$i=1,2,\ldots)$. Assume for contradiction that no $Y_i=X$.

\begin{lemma} There is a pre-compact open set $U$ in $X$ which is disjoint
from $Y_1$ and has non-empty intersection with infinitely many $Y_j$.
\end{lemma}

\begin{proof}

Fix $j>1$.  If  $\partial Y_j \subset Y_1$ then $X \setminus (Y_1
\cup Y_j)= \Int(Y_j) \cup (Y_j)^c$, which are disjoint open sets.
As intersections are small, this implies that $\Int(Y_j)$ is
empty, and hence $Y_j =\partial Y_j \subset Y_1$.  If this holds
for all $j>1$ then $X=Y_1$, contradiction.

Choose a $j>1$ for which $\partial Y_j$ is not contained in $Y_1$
and let $x$ be a point in $ \partial Y_j \setminus Y_1$.  Let $U$
be a pre-compact, connected neighborhood of $x$ which misses
$Y_1$.  If $U \cap Y_i$ is non-empty for infinity many $i$ then
the proof is complete.  Assume $U \subset Y_2 \cup Y_3 \cup \ldots
Y_n$. It is immediate that $U{\subset}\Int(Y_2{\cup} Y_3 {\cup}
\ldots Y_n)$ and that $Y_i{\cap}\Int(Y_2{\cup} Y_3 {\cup} \ldots
Y_n)=\Int(Y_i){\bigsqcup}\cup_{1<j<n,j{\neq}i} (Y_i \cap Y_j)$. It
follows that
$$U \setminus (\cup_{1< i<k \leq n} (Y_i \cap Y_k)){\subset}\bigsqcup_{1<i\leq n}(\Int(Y_i))$$
By Lemma \ref{lemma:trivial} the left hand side is connected,
hence contained in $\Int(Y_k)$ for some $k$. Thus implies $U
\subset Y_k$, which implies that $x \in \Int(Y_k)$ which is
disjoint from $Y_j$ and contradicts the choice of $x$.

\end{proof}

Let $U_1$ be the set produced by the lemma.  Apply the lemma to $U_1$ and
the cover $Y_2 \cup U_1, Y_3 \cup U_1, \ldots$ to get $U_2$ with
$\bar{U_2} \subset U_1$ and $U_2$ disjoint from $Y_2$.  Continue inductively
to get a sequence $U_i$ of pre-compact open sets with $\bar{U_i} \subset
U_{i-1}$ and $U_i$ disjoint from $Y_j$ for $j \leq i$.  The intersection
of the $U_i$ is equal to the intersection of the $\bar{U_i}$ and hence is
non-empty.  On the other hand, this intersection is disjoint for all the
$Y_i$, which is a contradiction.
\end{proof}

\begin{theorem} [Small Fixed Sets Implies Rigid]
\label{small} Let $\G$ be a countable group acting on $X$, a
connected, locally connected, locally compact Hausdorff
topological space. If $\Fix(\g)$ is small for all ${id} \neq \g
\in \G$ then $Inn(X,\G)=\G$, and every orbit equivalence is
equivariant.
\end{theorem}

These results have many applications.  For example:

\begin{corollary}
\label{corollary:hyperbolic}  Let $\G \subset SL_n(\Z)$ contain
only hyperbolic elements. The the orbit structure of $\G$ on
$\T^n$ is rigid.
\end{corollary}

As a simple consequence of Corollary \ref{corollary:hyperbolic},
for any standard hyperbolic action of ${\mathbb Z}^{n-1}$ on
${\mathbb T}^n$ we have $\Aut({\mathbb T}^n,{\mathbb
Z}^{n-1})=S_n{\ltimes}{\mathbb Z}^{n-1}$ where $S_n$ is the
symmetric group acting on ${\mathbb T}^n$ by permuting
coordinates.

\begin{corollary}  Let $\G \subset SL_n(\R)$ be any countable group.
The the orbit structure of $\G$ on $S^{n-1}$ is rigid.
\end{corollary}

Interestingly, this result is false if one allows orientation
reversing elements, see example \ref{example:circle}. That problem
does not occur over $\C$.

\begin{corollary}  Let $\G \subset PGL_n(\C)$ be any countable group.
The the orbit structure of $\G$ on $\C P^{n-1}$ is rigid.
\end{corollary}

There are also applications outside actions on manifolds.  For
example:

\begin{corollary} Let $\G$ be a word hyperbolic group which does not
split over any virtually cyclic group.  The orbit structure of
$\G$ on $\partial \G$ is rigid.  Furthermore if $\G$ and $\G'$ are
hyperbolic groups such that the boundary actions are orbit
equivalent then $\G$ and $\G'$ are commensurable.
\end{corollary}

\begin{proof}
Since the group does not split over a finite subgroup, the
boundary is connected by Stallings' Ends theorem \cite{St}.
Likewise, since the group does not split over any virtually
infinite cyclic group, the boundary has no cut points, by a
theorem of Bowditch \cite{Bow}.   Since the fixed set of any
element of $\G$ acting on the boundary is a pair of points, this
implies the fixed sets do not locally disconnect.  The boundary is
locally connected by \cite{BM}.    Thus the hypotheses of theorem
\ref{small} are satisfied, so $\Inn(\partial \G, \g) = \G$ and
every orbit equivalence is equivariant.   By \cite{P}, $\Out(\G)$
is finite, so $\G$ is finite index in $\Aut(\partial \G,\G)$.
\end{proof}

It would be interesting to know if similar results hold for
$CAT(0)$ groups.  The main obstacle seems to be determining when
the boundary of a flat can locally disconnect the boundary of the
group.  This might imply something like a splitting over an
abelian subgroup. However the situation is already unclear for
cutpoints, see e.g \cite{Sw}.

\section{Some Interesting Groups}
\label{groups}

At the boundary of rigidity, the groups of orbit equivalences can
be very interesting in their own right.  For example, Thompson's
group can be conjugated to a piecewise $PSL_2(\Z)$ action on the
circle (\cite{ghys}).    An interesting family of groups arise in
the case of $SL_n(\Z)$ acting on the $n$-torus.   Let $\G_n$ be
the group of homeomorphisms of $\T^n$ for which there is a
decomposition of $\T^n$ into finitely many pieces on which the
homeomorphism agrees with an element of $Gl_n(\Z)$.

\begin{lemma} $\G_n$ is the full group of Lipschitz orbit equivalences of the
$SL_n(\Z)$ action on $\T^n$.
\end{lemma}

\begin{proof}

 If $f$ is any Lipschitz orbit equivalence, then by the Lemma \ref{lemma:key}, $\T^n = \cup_{SL_n(\Z)} T_\g$,
with the union of the interiors dense.   Since $f$ is Lipschitz,
there are only finitely many elements of $SL_n(\Z)$ that can agree
with $f$ on an open dense set.  Thus there are $\g_1, \ldots,
\g_n$ in $SL_n(\Z)$ such that the union $T_{\g_i}$ is $\T^n$.
Since the intersections of two of these sets is contained in the
fixed set of some $\g_i {\g_j}^{-1}$, there is a finite union of
sub-tori such that $f$ restricted to any complementary component
is in $SL_n(\Z)$.

\end{proof}

This group is definitely bigger than $GL_n(\Z)$, see example \ref{adams}.

 Unlike Thompson's group, these groups are not simple.  There is a (split) surjection
from $\G_n$ to $GL_n(\Z)$ given by the action on $H^1(\T^n)$.
Explicitly, if $f$ is equal to $A_i$ on the set $X_i$, then the
image of $f$ in $GL_n(\Z)$ is $\Sigma A_i vol(X_i)$. Further, by
considering the germ of $f$ at $0$ one gets a homomorphism to
piecewise $Sl_n(\Z)$ acting on the $n-1$ sphere (by the double
cover of the standard projective action).

%\marginpar{Is this the right analogy with Thompson?  I thought
%Thompson was the kernel of some map from piecewise diadics to
%$\mathbb Z$. (I even think we figure this out when you visited me
%at Yale. Or rather, you figured it out and explained it to me.)}

\begin{lemma} Any subgroup of $\G_n$ is residually finite.
\end{lemma}

This follows from the interpretation as Lipschitz orbit
equivalences.  There are maps to finite permutation groups given
by restriction to any finite orbit.  Since finite orbits are
dense, the only thing in the kernel of all such maps is the
identity.

\begin{question}
\label{question:fg} Are the groups $\G_n$ finitely generated?
\end{question}

\begin{question} Is $\G_n$ the full group of orbit equivalences of $\T^n$ with the
standard $SL_n(\Z)$ action?
\end{question}

If the answer to Question \ref{question:fg} is no, one might try
to find interesting finitely generated subgroups of $\G_n$.

\medskip
\noindent
David Fisher\\
Department of Mathematics and Computer Science\\
Lehman College - CUNY\\
250 Bedford Park Boulevard W\\
Bronx, NY 10468\\
email:dfisher@lehman.cuny.edu\\
\\
\noindent
Kevin Whyte\\
Dept. of Mathematics\\
University of Illinois at Chicago\\
851 S. Morgan St.\\
Chicago, Il 60607\\
E-mail: kwhyte@math.uic.edu

\end{document}